\documentclass[12pt]{article}
\usepackage{amsfonts,amssymb,epsfig}
\usepackage{amsmath}

\date{}
\begin{document}
\newtheorem{df}{Definition}
\newtheorem{thm}{Theorem}
\newtheorem{lm}{Lemma}
\newtheorem{pr}{Proposition}
\newtheorem{co}{Corollary}
\newtheorem{re}{Remark}
\newtheorem{note}{Note}
\newtheorem{claim}{Claim}
\newtheorem{problem}{Problem}

\def\R{{\mathbb R}}

\def\E{\mathbb{E}}
\def\calF{{\cal F}}
\def\N{\mathbb{N}}
\def\calN{{\cal N}}
\def\calH{{\cal H}}
\def\n{\nu}
\def\a{\alpha}
\def\d{\delta}
\def\t{\theta}
\def\e{\varepsilon}
\def\t{\theta}
\def\pf{ \noindent {\bf Proof: \  }}
\def\trace{\rm trace}
\newcommand{\qed}{\hfill\vrule height6pt
width6pt depth0pt}
\def\endpf{\qed \medskip} \def\colon{{:}\;}
\setcounter{footnote}{0}

\def\Lip{{\rm Lip}}

\renewcommand{\qed}{\hfill\vrule height6pt  width6pt depth0pt}

\title{Approximate Gaussian isoperimetry for $k$ sets
}

\author{Gideon Schechtman\thanks{Supported by the
Israel Science Foundation.\ AMS subject classification: 60E15, 52A40
}}
\date{February 20, 2011}
\maketitle

\begin{abstract}

Given $2\le k\le n$, the minimal $(n-1)$-dimensional Gaussian measure of the union of the boundaries of $k$ disjoint sets of equal Gaussian measure in $\R^n$ whose union is $\R^n$ is of order $\sqrt {\log k}$.
A similar results holds also for partitions of the sphere $S^{n-1}$ into $k$ sets of equal Haar measure.

\end{abstract}

\section{Introduction}
Consider the canonical Gaussian measure on $\R^n$, $\gamma_n$. Given $k\in \N$ and $k$ disjoint measurable subsets of $\R^n$ each of $\gamma_n$ measure $1/k$ we can compute the $(n-1)$-dimensional Gaussian measure of the union of the boundaries of these $k$ sets. Below we shall make clear what exactly we mean by the $(n-1)$-dimensional Gaussian measure but in particular our normalization will be such that the $(n-1)$-dimensional Gaussian measure of a hyperplane at distance $t$ from the origin will be $e^{-t^2/2}$. The question we are interested in is what is the minimal value that this quantity can take when ranging over all such partitions of $\R^n$. As is well known, the Gaussian isoperimetric inequality (\cite{bo, st} implies that, for $k=2$, the answer is $1$ and is attained when the two sets are half spaces. The answer is also known for $k=3$ and $n\ge 2$ and is given by $3$ $2\pi/3$-sectors in $\R^2$ (product with $R^{n-2}$) (\cite{cch}). The value in question is then $3/2$.
If the $k$ sets are nice enough (for example if, with respect to the $(n-1)$-dimensional Gaussian measure, almost every point in the union of the boundaries of the $k$ sets belongs to the boundary of only two of the sets) then the quantity in question is bounded from below by $c\sqrt{log k}$ for some absolute $c>0$. This was pointed out to us by Elchanan Mossel. Indeed, by the Gaussian isoperimetric inequality, the boundary of each of the sets has measure at least $e^{-t^2/2}$ where $t$ is such that $\frac{1}{\sqrt{2\pi}}\int_t^\infty e^{-s^2/2}ds=1/k$. If $k$ is large enough $t$ satisfies
\[
\frac{e^{-t^2/2}}{\sqrt{2\pi}2t}<\frac{1}{k}<\frac{e^{-t^2/2}}{\sqrt{2\pi}t}
\]
which implies $\sqrt{\log k}\le t\le \sqrt{2\log k}$ and so the boundary of each of the $k$ sets has $(n-1)$-dimensional Gaussian measure at least $e^{-t^2/2}\ge \sqrt{2\pi}t/k\ge\sqrt{2\pi\log k}/k$. Under the assumption that the sets are nice we then get a lower bound of order $\sqrt{2\pi\log k}$ to the quantity we are after.

Of course the minimality of the boundary of each of the $k$ sets cannot occur simultaneously for even 3 of the $k$ sets (as the minimal configuration is a set bounded by an affine hyperplane) so it may come as a surprise that one can actually achieve a partition with that order of the size of the boundary. To show this is the main purpose of this note. It is natural to conjecture that, for $k-1\le n$ the minimal configuration is that given by the Voronoi cells of the $k$ vertices of a simplex centered at the origin of $\R^n$. So it would be nice to compute or at least estimate well what one gets in this situation. This seems an unpleasant computation to do. However, in Corollary \ref{co:1} below we compute such an estimate for a similar configuration - for even $k$ with $k/2\le n$, we look at the $k$ cells obtained as the Voronoi cells of $\pm e_i$, $i=1,\dots,k/2$ and show that the order of the $(n-1)$-dimensional Gaussian measure of the boundary is of order $\sqrt{\log k}$ and we deduce the main result of this note:

\medskip
\noindent {\bf Main Result} {\em Given even $k$ with $k\le 2n$, the minimal $(n-1)$-dimensional Gaussian measure of the union of the boundaries of $k$ disjoint sets of equal Gaussian measure in $\R^n$ whose union is $\R^n$ is of order $\sqrt {\log k}$.}
\medskip

In Corollary \ref{co:2} we deduce analogue estimates for the Haar measure on the sphere $S^{n-1}$.

This note benefitted from discussions with Elchanan Mossel and Robi Krauthgamer. I first began to think of the subject after Elchanan and I spent some time trying (alas in vain) to use symmetrization techniques to gain information on the (say, Gaussian) ``$k$-bubble" conjecture and some variant of it (see Conjecture 1.4 in \cite{im}). Robi asked me specifically the question that is solved here, with some possible applications to designing some algorithm in mind (but apparently the solution turned out to be no good for that purpose). I thank Elchanan and Robi also for several remarks on a draft of this note. I had also a third motivation to deal with this question. It is related to the computation of the dependence on $\e$ in (the probabilistic version of) Dvoretzky's theorem. It is too long to explain here, especially since it does not seem to lead to any specific result.

\section{Approximate isoperimetry for $k$ sets}

We begin with a simple inequality.
\begin{lm}\label {lm:1}
For all $\e>0$ if $C$ is large enough (depending on $\e$) then for all $k\in\N$
\[
\frac{1}{\sqrt{2\pi}}\int_{\sqrt{2\log \frac{k}{C}}
-1}^{\sqrt{2\log Ck}}\Big(\frac{1}{\sqrt{2\pi}}\int_{-s}^se^{-t^2/2}dt\Big)^{k-1}e^{-s^2/2}ds\ge\frac{1}{(2+\e)k}
\]
\end{lm}

\pf Let $g_1,g_2,\dots,g_k$ be independent identically distributed $N(0,1)$ variables. Then
\begin{equation}\label{eq:1}
\frac{1}{\sqrt{2\pi}}\int_0^\infty\Big(\frac{1}{\sqrt{2\pi}}\int_{-s}^se^{-t^2/2}dt\Big)^{k-1}e^{-s^2/2}ds=P(g_1\ge |g_2|,\dots,|g_k|)=\frac{1}{2k}.
\end{equation}
Also,
\begin{align}\label{eq:2}
\frac{1}{\sqrt{2\pi}}\int_0^{\sqrt{2\log \frac{k}{C}}
-1}\Big(\frac{1}{\sqrt{2\pi}}&\int_{-s}^se^{-t^2/2}dt\Big)^{k-1}e^{-s^2/2}ds\nonumber\\
&=\frac{1}{2k}\Big(\frac{2}{\sqrt{2\pi}}\int_0^se^{-t^2/2}\Big)^k\Big]_{s=0}^{\sqrt{2\log \frac{k}{C}}
-1}\\
&\le \frac{1}{2k}(1-\frac{2}{\sqrt{2\pi}}e^{-\log\frac{k}{C}})^k\le \frac{1}{2k}e^{-\frac{2C}{\sqrt{2\pi}}},\nonumber
\end{align}
and, for $C$ large enough,
\begin{align}\label{eq:3}
\frac{1}{\sqrt{2\pi}}\int_{\sqrt{2\log Ck} }^\infty&\Big(\frac{1}{\sqrt{2\pi}}\int_{-s}^se^{-t^2/2}dt\Big)^{k-1}e^{-s^2/2}ds\nonumber\\
&=\frac{1}{2k}\Big(\frac{2}{\sqrt{2\pi}}\int_0^se^{-t^2/2}\Big)^k\Big]_{s=\sqrt{2\log Ck}}^\infty\\
&\le \frac{1}{2k}\Big(1-\Big(1-\frac{2}{\sqrt{2\pi}}\int_{\sqrt{2\log Ck}}^\infty e^{-s^2/2}\Big)^k\Big)\le \frac{1}{2k}(1-e^{-1/C}).\nonumber
\end{align}
The Lemma now follows from (\ref{eq:1}),(\ref{eq:2}) and (\ref{eq:3}). \endpf

The next proposition is the main technical tool of this note. The statement involves the $(k-1)$-dimensional Gaussian measure of a certain subset of $\R^k$. We did not formally defined this notion for general sets yet (see Definition \ref{df:1} below) but the set we are talking about here is a subset of a hyperplane (through the origin of $\R^k$) and for such sets it just coincides with the canonical Gaussian measure, associated with this subspace, of the set in question.
\begin{pr}\label{pr:1}
For each $\e>0$ there is a $C$ such that for all $k\ge 2$, the $(k-1)$-dimensional Gaussian measure of the set $\{(t_1,t_2,\dots,t_k);\ t_1=t_2\ge |t_3|,\dots,|t_k|\}$ is bounded between $\frac{\sqrt{\pi\log \frac{k}{C}}-1}{(1+\e)2k(k-1)}$ and  $\frac{(1+\e)\sqrt{\pi\log Ck}}{2k(k-1)}$.
\end{pr}

\pf
The measure in question is
\[
\frac{1}{\sqrt{2\pi}}\int_0
^\infty\Big(\frac{1}{\sqrt{2\pi}}\int_{-s}^se^{-t^2/2}dt\Big)^{k-2}e^{-s^2}ds.
\]
Integration by parts (with parts $\Big(\frac{2}{\sqrt{2\pi}}\int_0^se^{-t^2/2}dt\Big)^{k-2}e^{-s^2/2}$ and $e^{-s^2/2}$) gives that this it is equal to
\begin{equation}\label{eq:4}
\frac{1}{2(k-1)}\int_0
^\infty\Big(\frac{2}{\sqrt{2\pi}}\int_0^se^{-t^2/2}dt\Big)^{k-1}se^{-s^2/2}ds.
\end{equation}
Now,

\begin{align}\label{eq:5}
\int_{\sqrt{2\log Ck}}
^\infty&\Big(\frac{2}{\sqrt{2\pi}}\int_0^se^{-t^2/2}dt\Big)^{k-1}se^{-s^2/2}ds\\
&=-\int_s
^\infty\Big(\frac{2}{\sqrt{2\pi}}\int_0^ue^{-t^2/2}dt\Big)^{k-1}e^{-u^2/2}du s\Big]_{s=\sqrt{2\log Ck}}^\infty\nonumber\\
&\phantom{==}+\int_{\sqrt{2\log Ck}}
^\infty\int_s^\infty\Big(\frac{2}{\sqrt{2\pi}}\int_0^ue^{-t^2/2}dt\Big)^{k-1}e^{-u^2/2}duds\nonumber\\
&\le\frac{\sqrt{2\pi}}{2k}(1-e^{-1/C})\sqrt{2\log Ck}+\int_{\sqrt{2\log Ck}}
^\infty\frac{\sqrt{2\pi}}{2k}(1-e^{-ke^{-s^2/2}})ds\label{eq:6},
\end{align}
where the estimate for the first term in (\ref{eq:6}) follows from (\ref{eq:3}) and of the second term follows from a similar computation to (\ref{eq:3}). Now (\ref{eq:6}) is at most
\begin{align}\label{eq:7}
\frac{\sqrt{2\pi}}{2Ck}\sqrt{2\log Ck}+\int_{\sqrt{2\log Ck}}
^\infty\frac{\sqrt{2\pi}}{2}e^{-s^2/2}ds\le \frac{\sqrt{2\pi}(\sqrt{2\log Ck}+1)}{2Ck}
\end{align}
and we conclude that
\begin{align}\label{eq:8}
\int_{\sqrt{2\log Ck}}
^\infty\Big(\frac{2}{\sqrt{2\pi}}\int_0^se^{-t^2/2}dt\Big)^{k-1}se^{-s^2/2}ds\le
\frac{\sqrt{2\pi}(\sqrt{2\log Ck}+1)}{2Ck}.
\end{align}
On the other hand
\begin{align}\label{eq:9}
&\int_0^{\sqrt{2\log Ck}}
\Big(\frac{2}{\sqrt{2\pi}}\int_0^se^{-t^2/2}dt\Big)^{k-1}se^{-s^2/2}ds\\
&\phantom{==}\le
\sqrt{2\log Ck}\int_0^\infty
\Big(\frac{2}{\sqrt{2\pi}}\int_0^s e^{-t^2/2}dt\Big)^{k-1}e^{-s^2/2}ds=\frac{\sqrt{2\pi}\sqrt{2\log Ck}}{2k}.\nonumber
\end{align}
Now, (\ref{eq:4}),(\ref{eq:8}) and (\ref{eq:9}) gives the required upper bound. The lower bound (which also follows from the Gaussian isoperimetric inequality) is easier. By Lemma \ref{lm:1}
\begin{align}\label{eq:10}
\frac{1}{2(k-1)}&\int_0
^\infty\Big(\frac{2}{\sqrt{2\pi}}\int_0^se^{-t^2/2}dt\Big)^{k-1}se^{-s^2/2}ds\\
&\ge \frac{1}{2(k-1)}\int_{\sqrt{2\log \frac{k}{C}}
-1}^{\sqrt{2\log Ck}}\Big(\frac{2}{\sqrt{2\pi}}\int_0^se^{-t^2/2}dt\Big)^{k-1}se^{-s^2/2}ds\\
&\ge \frac{{\sqrt{2\log \frac{k}{C}}-1}}{2(k-1)}\int_{\sqrt{2\log \frac{k}{C}}
-1}^{\sqrt{2\log Ck}}\Big(\frac{2}{\sqrt{2\pi}}\int_0^se^{-t^2/2}dt\Big)^{k-1}e^{-s^2/2}ds\\
&\ge \frac{{\sqrt{\pi\log \frac{k}{C}}-1}}{(1+\e)2k(k-1)}.
\end{align}
\endpf

To formulate Corollary \ref{co:1}, which is the main result here, in the most general setting we need to define the $(n-1)$-dimensional Gaussian measure of the boundary of a partition of $\R^n$ into $k$ sets.

\begin{df}\label{df:1}
Let $A_1,A_2,\dots,A_k$ be a partition of $\R^n$ into $k$ measurable sets. Put $A=\{A_1,A_2,\dots,A_k\}$ and denote
\[
\partial_\e A=\cup_{i=1}^k ((\cup_{j\not=i}A_j)_\e\setminus \cup_{j\not=i}A_j)
\]
(where $B_\e$ denotes the $\e$-neighborhood of the set $B$). We shall call $\partial_\e A$ the {\em $\e$-boundary} of $A$.
{\em The $(n-1)$-dimensional Gaussian measure of the boundary of $A$} will be defined and denoted by
\[
\gamma_{n-1}(\partial A)=\liminf_{\e\to 0}\frac{\gamma_n(\partial_\e A)-\gamma_n(A)}{\sqrt{2/\pi}\e}.
\]
\end{df}
Note that we do not define the boundary of the partition, only the measure of the boundary. We would like however that in simple cases, when the boundary and its $(n-1)$-dimensional Gaussian measure are well understood, this definition will coincide with the classical one. In particular notice that if the partition is into two sets which are separated by a hyperplane at distance $t$ from the origin the definition says that the $(n-1)$-dimensional Gaussian measure of the boundary is $e^{-t^2/2}$ and in particular when $t=0$ the measure is 1 which coincides with what we understand as the classical $\gamma_{n-1}$ measure of a hyperplane through 0. This is why the factor $\sqrt{2/\pi}$ is present in the definition above.

\begin{co}\label{co:1} For some universal constants $0<c<C<\infty$ and all $k=2,3,\dots$,\\
(1) If $A=\{A_1,A_2,\dots,A_k\}$ is a partition of $\R^n$ into $k$ measurable sets each of $\gamma_n$ measure $1/k$. Then $\gamma_{n-1}(\partial A)\ge c \sqrt{\log k}$.\\
(2) If $k\le n$, there is a partition $A=\{A_1,A_2,\dots,A_{2k}\}$ of $\R^n$ into $2k$ measurable sets each of $\gamma_n$ measure $1/2k$ such that $\gamma_{n-1}(\partial A)\le C \sqrt{\log k}$.
\end{co}

(1) follows very similarly to the argument in the introduction, except that there is no need for the boundary to be nice anymore: By the Gaussian isoperimetric inequality, for each $\e>0$ and each $i=1,\dots,k$,
\[
\gamma_n((\cup_{j\not=i}A_j)_\e\setminus \cup_{j\not=i}A_j)\ge\frac{1}{\sqrt{2\pi}}\int_t^{t+\e}e^{-s^2/2}ds,
 \]
 where $t$ is such that $\frac{1}{\sqrt{2\pi}}\int_t^\infty e^{-s^2/2}ds=1/k$. If $\e$ is small enough, the argument in the introduction gives that the integral in question is of order $\e\frac{\sqrt{\log k}}{k}$. Since the $k$ sets $(\cup_{j\not=i}A_j)_\e\setminus \cup_{j\not=i}A_j$ are disjoint, we deduce (1).
(2) follows directly from Proposition \ref{pr:1} since the boundary of the partition into the Voronoi cells corresponding to $\{\pm e_i\}_{i=1}^k $ is contained in the union of $k(k-1)$ hyperplans through zero and thus $\gamma_{n-1}(\partial A)$ coincide with the classical $\gamma_{n-1}(\partial A)$ which is what is estimated in Proposition \ref{pr:1}.

A similar result to Corollary \ref{co:1} holds on the $n$-dimensional sphere, $S^{n-1}$ with its normalized Haar measure $\sigma_n$. One defines the $\e$-boundary of a partition $A$ of the sphere in a similar way to the first part of Definition \ref{df:1} (using, say, the geodesic distance to define the $\e$-neighborhood of a set). Then one defines the {\em $(n-1)$-dimensional Haar measure of the boundary of $A$} by
\[
\sigma_{n-1}(\partial A)=\liminf_{\e\to 0}\frac{\sigma_n(\partial_\e A)-\sigma_n(A)}{\sqrt{2n/\pi}\e}.
\]
The choice of the normalization constant $\sqrt{2n/\pi}$ was made so that if the partition is into two sets separated by a hyperplane then the measure of the boundary (which ``is" $S^{n-2}$) will be 1. The proof can be obtained from that of Corollary \ref{co:1} by a standard reduction, using the fact that if $(g_1,\dots,g_n)$ is a standard Gaussian vector then the distribution of  $(\sum g_i^2)^{-1}(g_1,\dots,g_n)$ is $\sigma_n$.

\begin{co}\label{co:2} For some universal constants $0<c<C<\infty$ and all $k=2,3,\dots$,\\
(1) If $A=\{A_1,A_2,\dots,A_k\}$ is a partition of $S^{n-1}$ into $k$ measurable sets each of $\sigma_n$ measure $1/k$. Then $\sigma_{n-1}(\partial A)\ge c \sqrt{\log k}$.\\
(2) If $k\le n$, there is a partition $A=\{A_1,A_2,\dots,A_{2k}\}$ of $S^{n-1}$ into $2k$ measurable sets each of $\sigma_n$ measure $1/2k$ such that $\sigma_{n-1}(\partial A)\le C \sqrt{\log k}$.
\end{co}

\begin{re}
It may be interesting to investigate what happens when $k>>n$. In particular, if $k=2^n$ then the  partition of $R^n$ into its $k=2^n$ quadrants satisfy that the $\gamma_{n-1}$ measure of its boundary (consisting of the coordinates hyperplanes) is $n=\log k$. Is that the best (order) that can be achieved?
\end{re}


\noindent Gideon Schechtman\newline Department of
Mathematics\newline Weizmann Institute of Science\newline Rehovot,
Israel\newline E-mail: gideon.schechtman@weizmann.ac.il

\end{document}